\documentclass[12pt]{article}
\usepackage[english]{babel}
\usepackage{amsmath,amsthm}
\usepackage{amsfonts}

\theoremstyle{definition}

\theoremstyle{remark}

\numberwithin{equation}{section}

\begin{document}

\begin{center}
\bfseries  Hardy-Littlewood Conjecture and Exceptional real Zero  \\
\end{center}

\begin{center}

\vspace{5mm}
JinHua  Fei\\
\vspace{5mm}

ChangLing  Company of Electronic Technology    \, Baoji \,  Shannxi  \,  P.R.China \\
\vspace{5mm}

E-mail: feijinhuayoujian@msn.com \\

\end{center}

\vspace{3mm}

{\bfseries Abstract.}\, In this paper, we assume that  Hardy-Littlewood Conjecture, we got a better upper bound of the exceptional real zero for a class of module. \\


{\bfseries MR(2000)  Subject  Classification \quad 11P32,  11M20 } \\

\vspace{8mm}

In this paper, we reveal the relationship between  Hardy-Littlewood Conjecture and the exceptional real zero. Assume that Hardy-Littlewood Conjecture, we obtained a good results of the exceptional real zero. \\

In this paper, I generalize the results of my paper "An Application of Hardy-Littlewood Conjecture". $q$ is the  prime number and $ q \equiv 3 \,(mod \,4)$ is improved to $q$ is  odd square-free and $ q \equiv 3 \,(mod \,4).$ we know that the module $q$ of the exceptional primitive real character are square-free when q is odd positive integer.\\

First, we give Hardy-Littlewood Conjecture.\\

{\bfseries\footnotesize Hardy-Littlewood Conjecture.}  When $N$ is even integer and $N\geq 6$, we have

$$ \sum_{3\leq p_1,p_2\leq N \atop p_1 +p_2 =N} 1  \approx  \frac{N}{\varphi (N)} \prod_{p \, \dagger N} \left( 1-\frac{1}{(p-1)^2}   \right ) \frac{N}{\log^2 N } $$ \\

where $p_1,p_2,p$ are the prime numbers, $\varphi(n) $ is Euler function.\\

Under the above conjecture, we have the following theorem \\

{\bfseries\footnotesize Theorem.} Let $q$ is odd square-free and $q \equiv \,3 \,(mod\, 4),$ it has exceptional real character $\chi$,  and  its Dirichlet $L(s,\chi)$  function has an exceptional real zero $\beta$.  If Hardy-Littlewood Conjecture is correct, then there is a  positive constant $c$ , we have\\

$$\beta \leq 1- \frac{c}{\log^2 q}$$ \\

Now, we do some preparation work.\\

{\bfseries\footnotesize Lemma 1 .} Assuming that the Hardy-Littlewood Conjecture. Let $ n $ is any positive integer, $ q $ is odd  integer and large sufficiently, then\\

$$ \sum_{3\leq p_1,p_2\leq 2nq \atop p_1 +p_2 =2nq } 1  \geq  \frac{q}{\varphi (q)}  \frac{2\,n\,q \,d}{(\log 2nq)^2 } $$ \\

where

$$ d = \prod_{3\leq p } \left( 1-\frac{1}{(p-1)^2}   \right ) $$  \\

and $ p_1,p_2, p   $ are the prime numbers.\\

{\bfseries\footnotesize  Proof .} By Hardy-Littlewood Conjecture, when $N$ is large sufficiently, we have

$$ \sum_{3\leq p_1,p_2\leq N \atop p_1 +p_2 =N} 1  \geq   \frac{N}{2 \varphi (N)} \prod_{p \, \dagger N} \left( 1-\frac{1}{(p-1)^2}   \right ) \frac{N}{\log^2 N } $$ \\

because

$$\prod_{p \, \dagger N} \left( 1-\frac{1}{(p-1)^2}   \right ) \geq \prod_{3 \leq p \,} \left( 1-\frac{1}{(p-1)^2}   \right )  $$ \\

and

$$  \frac{2nq}{\varphi (2nq)} = \prod_{p|2nq} \frac{p}{p-1 }  \geq \prod_{p|2q} \frac{p}{p-1 } = \frac{2q}{\varphi (q)} $$\\

We choose $ N=2nq  $ , This completes the proof of Lemma 1 .  \\

{\bfseries\footnotesize Lemma 2.} Let $ m $  is  positive integer and $ n $ is integer, then\\

\[  \sum_{k=1}^m e\left(\frac{kn}{m}\right) = \left\{
  \begin{array}{ll}
  m  \quad if \,\, n \equiv 0 \,\, (mod \,\, m)\\
  0  \quad  otherwise
  \end{array}
  \right.
\]\\

where $e(x)=e^{2\pi i x}$.  The lemma 2 is obvious\\

{\bfseries\footnotesize Lemma 3.}  Let $c_1$ be the positive constant. if $ (a,q)=1 $, then

$$ \pi (x;q,a) =  \frac{Li x}{\varphi (q)} - \frac{\chi(a)}{\varphi (q)} \int_2^x \frac{u^{\beta - 1}}{\log u} du +  O \left(x \exp(-c_1\sqrt{\log x})\right)$$ \\

when there is an exceptional character $ \chi $ modulo $ q $  and $ \beta $ is the concomitant zero. Where $ Li x =\int_2^x \frac{du}{\log u} $ and $ \exp (x) = e^x$\\

 The lemma 3 follows from the  References [2],   Corollary 11.20 of the  page  381   \\

It is easy to see that

$$ Li x =\int_2^x \frac{du}{\log u} = \frac{x}{\log x}  + O\left( \frac{x}{\log^2 x}   \right) $$  \\

and

$$ \int_2^x \frac{u^{\beta - 1}}{\log u} du  = \frac{x^\beta}{\beta \log x} + O\left(  \frac{x^\beta}{\log^2 x }  \right) $$  \\

{\bfseries\footnotesize Lemma 4.}  if $ \chi $ is a primitive character modulo $ m, $   then

$$ \sum_{k=1}^m \chi(k) e\left(\frac{nk}{m}\right) = \overline{\chi} (n) \tau (\chi) $$

where $ \tau (\chi) = \sum_{k=1}^m \chi(k) e(\frac{k}{m}).   $\\

The lemma 4 follows from the  References [1], the page 47. \\

{\bfseries\footnotesize Lemma 5.}   if $ m $ is odd square-free and $\chi$ is a primitive real character modulo $ m, $  then

\[  \tau (\chi) = \left\{
  \begin{array}{ll}
  \sqrt{m}  \quad  \, if \,\, m \equiv 1 \,\, (mod \,\, 4)\\
  i \sqrt{m}  \quad  if \,\, m \equiv 3 \,\, (mod \,\, 4)
    \end{array}
  \right.
\]

The lemma 5 follows from the  References [1], the theorem 3.3 of the page 49.      \\

{\bfseries\footnotesize Lemma 6.}  We give the value of two sums, they are used in the proof of the Theorem.  \\

(1) \\

$$  \sum_{k=1}^q \left( \sum_{a=1 \atop (a,q)=1}^{q-1} e \left ( \frac{ak}{q} \right )  \right)^2 = \sum_{k=1}^q \sum_{a=1 \atop (a,q)=1}^{q-1} e \left ( \frac{ak}{q} \right )  \sum_{b=1 \atop (b,q)=1}^{q-1} e \left ( \frac{bk}{q} \right )  $$ \\

$$ = \sum_{k=1}^q \sum_{a=1 \atop (a,q)=1}^{q-1}   \sum_{b=1 \atop (b,q)=1}^{q-1} e \left ( \frac{(a+b)k}{q} \right ) =   \sum_{a=1 \atop (a,q)=1}^{q-1}   \sum_{b=1 \atop (b,q)=1}^{q-1}  \sum_{k=1}^q e \left ( \frac{(a+b)k}{q} \right )$$ \\

$$ =  q \sum_{a=1 \atop (a,q)=1}^{q-1}   \sum_{b=1 \atop (b,q)=1,\, a+b=q  }^{q-1} 1 = q \, \varphi (q)  \quad \quad \quad \quad \quad \quad \quad \quad \quad \quad  $$     \\

(2)\\

$$  \sum_{k=1}^q \chi(k) \sum_{a=1 \atop (a,q)=1}^{q-1} e \left( \frac{ak}{q} \right) = \sum_{a=1 \atop (a,q)=1}^{q-1}  \sum_{k=1}^q \chi(k) e \left( \frac{ak}{q} \right) = \tau(\chi) \sum_{a=1 \atop (a,q)=1}^{q-1} \overline{\chi}(a) = 0  $$\\

where $\chi $ is the primitive character modulo $ q, $ \\

This completes the proof of Lemma 6 .  \\

{\bfseries\footnotesize PROOF OF THEOREM.}  \\

 \emph{The first part}.\\

 By Lemma 2, when $ x \geq q^4, $  we have

$$  \sum_{k=1}^q \left(  \sum_{3\leq p \leq x} e\left(\frac{kp}{q}\right)     \right)^2 =   \sum_{k=1}^q   \sum_{3\leq p_1 \leq x} \sum_{3\leq p_2 \leq x} e\left(\frac{k(p_1+ p_2)}{q}\right)  $$  \\

$$ =      \sum_{3\leq p_1 \leq x} \sum_{3\leq p_2 \leq x} \sum_{k=1}^q  e\left(\frac{k(p_1+ p_2)}{q}\right) = q  \sum_{3\leq p_1,  p_2 \leq x \atop p_1 + p_2 \equiv \, 0 \,(q)} 1  \geq q \sum_{n=1}^{[\frac{x}{2q}]} \sum_{3\leq p_1,  p_2 \leq x \atop p_1 + p_2 = 2nq } 1  $$  \\

by Lemma 1 , the above formula

$$  \geq \frac{q^2}{\varphi(q)} \sum_{n=1}^{[\frac{x}{2q}]}  \frac{\ 2nq d}{\log^2 2nq} \geq \frac{q^2}{\varphi(q)}  \sum_{n=1}^{[\frac{x}{2q}]}  \frac{2nqd}{\log^2 x }  \geq   \frac{2 d q^3}{ \varphi (q)  \log^2 x}  \sum_{n=1}^{[\frac{x}{2q}]} n $$  \\

$$  = \frac{2 d q^3 }{ \varphi (q)  \log^2 x} \cdot \,\frac{[\frac{x}{2q}]([\frac{x}{2q}]+1 )}{2}  \geq \frac{q \,d \,x^2}{ 4 \varphi (q) \log^2 x} + O\left(\frac{xq^2}{ \varphi (q) \log^2 x}\right)  $$  \\

\emph{The second part.}\\

When $ 1\leq k \leq q ,$  we have \\

$$ \sum_{3 \leq p \leq x}  e\left(\frac{pk}{q}\right) = \sum_{3 \leq p \leq x \atop (p,q) =1  }  e\left(\frac{pk}{q}\right) + O(\log q) = \sum_{a=1 \atop (a,q)=1 }^q  e\left(\frac{ak}{q}\right) \sum_{3\leq p \leq x \atop p \equiv a (q)} 1  + O(\log q) $$  \\

by Lemma 3 and Lemma 4,  the above formula

$$ = \sum_{a=1 \atop (a,q)=1 }^q  e\left(\frac{ak}{q}\right) \left( \frac{Li x}{\varphi (q)} - \frac{\chi(a)}{\varphi (q)} \int_2^x \frac{u^{\beta - 1}}{\log u} du +  O \left(x \exp(-c_1\sqrt{\log x})\right)    \right)  + O(\log q)$$ \\

$$ =   \frac{  Li x}{ \varphi (q)}  \sum_{a=1 \atop (a,q)=1 }^q  e\left(\frac{ak}{q}\right)   - \frac{ \tau (\chi) \chi(k)}{\varphi (q )} \int_2^x \frac{u^{\beta - 1}}{\log u} du +  O \left( q x \exp(-c_1\sqrt{\log x})\right)$$  \\

where $\chi$ is the exceptional primitive real character modulo $ q. $\\

therefore\\

$$   \left( \sum_{3 \leq p \leq x}  e\left(\frac{pk}{q}\right) \right)^2   =  \left( \frac{  Li x}{ \varphi (q)}  \sum_{a=1 \atop (a,q)=1 }^q  e\left(\frac{ak}{q}\right)  \right)^2  $$  \\

$$ - \frac{2 \tau(\chi) \chi(k)Lix }{\varphi^2 (q)}  \left( \sum_{a=1 \atop (a,q)=1 }^q  e\left(\frac{ak}{q}\right)\right)\int_2^x \frac{u^{\beta - 1}}{\log u} du  +\left( \frac{\chi(k)\tau (\chi)}{\varphi(q)} \int_2^x \frac{u^{\beta - 1}}{\log u} du \right)^2  $$\\

 $$ + O \left( q^2 x^2 \exp(-c_2\sqrt{\log x}     \right)  \qquad \qquad \qquad \quad \quad \quad   \quad \quad \quad \quad \quad \quad $$  \\

By Lemma 5 and Lemma 6, we have \\

 $$ \sum_{k=1}^{q} \left( \sum_{3 \leq p \leq x}  e\left(\frac{pk}{q}\right) \right)^2 = \frac{q}{\varphi(q)}\left(Li x \right)^2  - \frac{q}{\varphi(q)} \left(\int_2^x \frac{u^{\beta - 1}}{\log u} du \right)^2  $$  \\

$$  +  O \left(   q^3 x^2 \exp(-c_2\sqrt{\log x}     \right) \qquad  \qquad  \qquad   \qquad   \qquad  \qquad   \qquad $$ \\

$$ =\frac{q \,x^2}{\varphi(q) \log^2 x} -  \frac{q\, x^{2\beta}}{ \varphi(q) \beta^2  \log^2 x} +  O \left( \frac{q \,x^2}{\varphi(q) \log^3 x}  + q^3 x^2 \exp(-c_2\sqrt{\log x}     \right)  $$  \\

We synthesize the first part and second part, we have

$$ \frac{q \,d \,x^2}{ 4 \varphi (q) \log^2 x}  \leq \frac{q \,x^2}{\varphi(q) \log^2 x}  - \frac{q \,x^{2\beta}}{\varphi (q)\beta^2 \log^2 x }   $$  \\

$$ + O \left(  \frac{x^2}{ \log^3 x  } +  \frac{xq^2}{ \varphi (q) \log^2 x} +  q^3 x^2 \exp(-c_2\sqrt{\log x})\right)  $$ \\

$$ \frac{\,d \,x^2}{ 4  \log^2 x}  \leq \frac{\,x^2}{ \log^2 x}  - \frac{ \,x^{2\beta}}{\beta^2 \log^2 x }  + O \left(  \frac{x^2}{ \log^3 x  } +  \frac{x\,q}{ \log^2 x} +  q^3 x^2 \exp(-c_2\sqrt{\log x})\right)   $$  \\

$$ \frac{\,d }{ 4 }  \leq 1  - \frac{ \,x^{2\beta-2}}{\beta^2 }  + O \left(  \frac{1}{ \log x  } +  \frac{q}{ x} +  q^3 \exp(-c_3\sqrt{\log x})\right)   $$  \\

we take $ \log x = ( \frac{4}{c_3} \log q )^2 , $ then

$$ x^{2\beta-2}  \leq  1-\frac{d }{4}+ \frac{c_4}{\log^2 q} $$  \\

we take $ \log q \geq  \sqrt{\frac{8 c_4}{d}}, $  then

$$ x^{2\beta-2}  \leq  1-\frac{d}{8} $$  \\

$$  \beta-1 \leq  \frac{\log ( 1-\frac{ d }{8}) }{2\log x} = -  \frac{\log (\frac{8}{8-d} ) }{2\log x} $$ \\

therefore

$$   \beta \leq 1- \frac{c}{\log^2 q}  $$  \\

This completes the proof of  Theorem. \\

Because $ \tau(\chi_4) = 2i  $ and $ \tau(\chi_8) = 2\sqrt{2} $ , when $q$ is  odd square-free, by  same method, for module $4q$ ,$q\equiv 1 (mod 4)$ and  module $8q$ ,$q\equiv 3 (mod 4)$ , we have the same conclusion.

\vspace{20mm} \centerline{ REFERENCES } \vspace{5mm}

[1]  Henryk Iwaniec, Emmanuel Kowalski, { \itshape Analytic Number Theory,} American mathematical Society, 2004. \\

[2]  Hugh L. Montgomery,  Robert C. Vaughan, {\itshape Multiplicative Number Theory I. Classical Theory, } Cambridge University Press, 2006.  \\

\end{document}